\documentclass[11pt]{article}
\usepackage{amsthm, amsfonts, url}
\usepackage{fullpage}
\usepackage{amscd,amsmath,amssymb}

\newtheorem{definition}{Definition}
\newtheorem{theorem}{Theorem}
\newtheorem{lemma}{Lemma}

\newtheorem{corollary}{Corollary}

\def\al{\alpha}

\def\La{\Lambda}
\def\la{\lambda}

\def\kappa{\varkappa}

\def\C{{\mathbb C}}
\def\y{{\mathbb Y}}
\def\Z{{\mathbb Z}}
\def\yn{{\mathbb Y}_n}
\def\sn{{\mathfrak S}_n}
\def\si{{\mathfrak S}_{\mathbb N}}
\def\gn{{\mathbb C}[\sn]}

\def\gzn{\operatorname{GZ}_n}

\def\N{{\mathbb N}}

\def\Prob{\operatorname{Prob}}

\def\be{\begin{equation}}
\def\ee{\end{equation}}

\def\Fnk{F_{n,k}}

\urldef\vershik\url{vershik@pdmi.ras.ru}
\urldef\natalia\url{natalia@pdmi.ras.ru}

\author {A.~M.Vershik\thanks{%
St.~Petersburg Department of Steklov Institute of Mathematics.
E-mail: \vershik, \natalia. Supported by the
CRDF grant  RUM1-2622-ST-04 and INTAS grant 03-51-5018.%
} \and N.~V.~Tsilevich\footnotemark[1]}
\title {Markov measures on Young tableaux
and induced representations of the infinite symmetric group}

\date{}

\begin{document}
\maketitle

\begin{abstract}
We show that the class of so-called Markov representations of
the infinite symmetric group $\si$, associated with Markov
measures on the space of infinite Young tableaux, coincides 
with the class of simple representations, i.e., inductive
limits of representations with simple spectrum.
The spectral measure of an arbitrary representation of $\si$
with simple spectrum is equivalent to a multi-Markov measure 
on the space of Young tableaux.
We also show that the representations of $\si$ induced from the
identity representations of two-block
Young subgroups are Markov and find explicit formulas
for the transition probabilities of the corresponding Markov measures.
The induced representations are studied with the help of the tensor
model of two-row representations of the symmetric groups;
in particular, we deduce explicit formulas for the Gelfand--Tsetlin basis
in the tensor models.We show that the class of so-called Markov representations of
the infinite symmetric group $\si$, associated with Markov
measures on the space of infinite Young tableaux, coincides 
with the class of simple representations, i.e., inductive
limits of representations with simple spectrum.
The spectral measure of an arbitrary representation of $\si$
with simple spectrum is equivalent to a multi-Markov measure 
on the space of Young tableaux.
We also show that the representations of $\si$ induced from the
identity representations of two-block
Young subgroups are Markov and find explicit formulas
for the transition probabilities of the corresponding Markov measures.
The induced representations are studied with the help of the tensor
model of two-row representations of the symmetric groups;
in particular, we deduce explicit formulas for the Gelfand--Tsetlin basis
in the tensor models.
\end{abstract}

\section{Introduction}
Consider the Young graph $\mathbb Y$, i.e., the $\mathbb
Z_+$-graded graph of Young diagrams: the vertices of its 
$n$th level $\yn$ are Young diagrams with
$n$ cells, and edges connect two vertices of neighboring levels
whenever the corresponding diagrams differ by one cell.
The space $T$ of Young tableaux, i.e., infinite paths in the Young graph, 
is a totally disconnected (nonstationary) Markov compactum, 
and Markov measures on $T$ are defined in the usual way. In other words,
for every fixed diagram $\lambda_n\in\yn$, the conditional 
measure on the set of tableaux ``passing through'' this diagram
at level $n$ is the direct product of the conditional measures
on the ``past'' (before the moment $n$) and ``future'' (after the moment $n$), 
i.e., the past and the future of a path are independent under a fixed 
``present'' diagram. Such a measure is determined by the conditional
transition (or cotransition) probabilities. This in turn means 
that the sequence of random variables
$X_n$ (where $X_n(t)$ is the content of the $n$th cell of a tableau
$t$) has a special structure:
the distribution of  $X_{n+1}$ depends not on the all previous values
$X_k$, but only on the ``border'' ones (i.e., only on the diagram).

The tail partition $\xi$ on the space of paths $T$ is the partition
into classes of confinal paths. A measure 
$\mu$ (not necessarily Markov) on $T$ is called ergodic
(with respect to the tail partition) if every measurable $\xi$-set
has either zero or full $\mu$-measure. A measure
$\mu$ is called quasi-invariant if it is quasi-invariant in the ordinary sense
under any transformations changing the beginnings of tableaux.
In other words, for almost every tableau
$t$ and every $n\in\N$, the conditional measure on the space of tableaux
coinciding with $t$ starting from the moment $n$ is strictly positive. 

Recall (see \cite{James, JK, Itogi}) that the infinite symmetric group
$\si$ is the inductive limit of the finite symmetric groups $\sn$
with the natural embeddings; the Young diagrams with $n$ cells
parametrize the irreducible complex representations of $\sn$;
the Young graph is the branching graph of irreducible
representations of the symmetric groups (= Bratteli diagram
of the group algebra of $\si$).
It is well known (see \cite{Itogi}) that each ergodic quasi-invariant
measure on the space of infinite tableaux $T$ and each $1$-cocycle
(on the tail equivalence relation in this space) with values
in the group of complex numbers of modulus $1$ determine 
an irreducible representation of
$\si$ with simple spectrum (with respect to the Gelfand--Tsetlin algebra)
by the formula~(\ref{repr}) (see \S\ref{sec:markov}).
Conversely, each irreducible representation with simple spectrum
can be realized in this form.
Two such representations are equivalent if and only if
the corresponding measures are mutually absolutely continuous and the cocycles 
are cohomological. Not all irreducible representations
can be obtained in this way, there exist irreducible representations
with nonsimple spectrum; but the representations in question form a natural
and important class of irreducible representations, which is defined in
invariant terms, i.e., closed under all automorphisms of the group $\si$.

A representation of the infinite symmetric group  $\si$ determined by a Markov measure is called
a {\it Markov representation}; thus a Markov representation can be determined 
by the transition probabilities of adding a new cell to a Young diagram 
(or the cotransition probabilities of deleting a cell).

We define a {\it generalized Markov measure} on the space $T$
as a measure for which there exists an increasing sequence
$N_1,N_2,{\ldots}$ of positive integers such that for every $n$
the conditional measures on the past and on the future are independent
if one fixes diagrams at the interval 
$N_n,{\ldots} ,N_{n+1}-1$.

By analogy with the corresponding theorem from 
\cite{V} on quasi-invariant measures on the spaces of sequences,
one can use the martingale convergence theorem
to prove the
following result: every quasi-invariant ergodic measure is
mutually absolutely continuous with a generalized Markov
quasi-invariant ergodic measure. Thus, for the study of representations
with simple spectrum, it suffices to consider only quasi-invariant
generalized Markov measures.

Let us say that a representation of $\si$ is {\it simple} if it
is the inductive limit of a sequence of representations of the finite
symmetric groups
$\sn$, $n=1,2,{\ldots}$, each having a simple spectrum, i.e.,
having no multiplicities in the decomposition
into irreducible representations. Simple representations have a simple spectrum,
but they do not exhaust all representations with simple spectrum.
The class of simple representations  essentially
depends on the approximation of the group $\si$
by finite groups (see the end of Sec.~\ref{sec:markov});
but in this paper we consider only the standard approximation, so that
we use the term ``simple representation'' omitting the reference to the
approximation. It turns out  (Theorem~\ref{th:mark})
that the class of simple representations coincides with the class of
Markov representations. This simple yet important fact relates
the notions that are quite different at first sight.

In this paper, we consider only two-row Young diagrams, which
fill only a small part of the Young graph, namely, the ``half''
of the Pascal graph
(the vertices are indexed by pairs of positive integers
$(n,k)$, where $n$ is the level number, $k$ is the vertex number, 
and $k \le [n/2]$). A path in this graph 
is a trajectory of a random walk on the semilattice
${\mathbb Z}_+$ with reflecting barrier at the origin.
All previous definitions (Markov measures, representations, etc.)
can be restricted to this case.

Representations of the finite symmetric groups determined by two-row
diagrams are widespread 
in applications; they have a convenient model --- the so-called
tensor realization, where the representation space is a space
of finite-dimensional tensors; the rank of tensors
does dot exceed half the dimension
(see \cite{Nik}). We use this model below.

Our main result (Theorem~\ref{th:spectral}) claims that
the representations of the infinite symmetric group induced from
two-block
Young subgroups are simple, and, consequently, Markov;
we also find explicit formulas for the corresponding
Markov spectral measures, which are
the laws of remarkable and natural random walks on the semilattice.

This result illustrates the thesis from
\cite{Unlim} on the contemporary fruitful rapprochement between
algebraic and probabilistic constructions: the notion of induced
representation is one of the central notions in representation theory,
while that of Markov measures plays an important role in probability theory.
In our construction, both notions are amazingly intertwined.

\medskip
The paper is organized as follows.
In Sec.~\ref{sec:not}, we recall basis notions and
introduce notation related to
the Young graph and representation theory of the symmetric groups,
in particular, the notion of Markov measure on the space of Young tableaux,
which plays a key role in our considerations.
In Sec.~\ref{sec:markov}, we introduce the notions of Markov vector
and Markov representation and prove that
the classes of Markov and simple representations coincide.
Section~\ref{sec:tensor} describes the so-called tensor model
of two-row representations of the finite symmetric groups. In Sec.~\ref{sec:gz},
we find explicit formulas for the Gelfand--Tsetlin basis in the tensor model.
Finally, in Sec.~\ref{sec:induced} we prove that the representations
of the infinite symmetric group induced from the identity representations
of two-block subgroups are Markov and
explicitly find the corresponding spectral measures.

\section{The Young graph and Markov measures on the space of its paths}
\label{sec:not}

In this section, we recall necessary notions related to
the Young graph and the representation theory of the symmetric groups
(see, e.g., \cite{James, JK} and \cite{GB}); in particular,
the notion of a Markov measure on the space of Young tableaux.

We denote by $\sn$ the symmetric group of degree
$n$ and by
$\gn$ the group algebra of $\sn$.

The irreducible representations of the symmetric group
$\sn$ are indexed by the set
$\yn$ of Young diagrams with $n$ cells. Let
$\pi_\la$ be the irreducible unitary representation of
$\sn$ corresponding to a diagram
$\la\in\yn$, and let $\dim\la$ be 
the dimension of $\pi_\la$.

The branching of irreducible representations of the symmetric groups
is described by the Young graph $\y$. The set of vertices of the $\Z_+$-graded 
graph $\y$ is
$\cup_n\yn$, and two vertices $\mu\in{\mathbb Y}_{n-1}$ and $\la\in\yn$
are joined by an edge if and only if $\mu\subset\la$.
By definition, the zero level ${\mathbb Y}_{0}$ consists of
the empty diagram $\emptyset$.

Denote by $T_\la$ the set (consisting of $\dim\la$ elements)
of Young tableaux of shape
$\la\in\yn$, or, which is the same, the set of paths in the Young graph
from the empty diagram $\emptyset$ to $\la$.
Let
$T_n=\cup_{\la\in\yn}T_\la$
be the set of Young tableaux with $n$ cells.

According to the branching rule for irreducible representations
of the symmetric groups, the space $V_\la$ of the irreducible
representation  $\pi_\la$ decomposes into the sum of one-dimensional subspaces
indexed by the tableaux $u\in T_\la$. The basis $\{h_u\}_{u\in T_\la}$
consisting of vectors of these subspaces is called the
{\it Gelfand--Tsetlin basis}. It is an eigenbasis for the {\it Gelfand--Tsetlin
algebra} $\gzn$, the subalgebra in $\gn$ generated by the centers
$Z[{\mathfrak S}_1], Z[{\mathfrak S}_2],{\ldots} ,Z[{\mathfrak S}_n]$
(see \cite{VO}). 

Now let
$\si=\cup_{n=1}^\infty\sn=\varinjlim\sn$ be the infinite symmetric group
with the fixed structure of an inductive limit of finite groups.

Denote by $T=\varprojlim T_n$ the space of infinite Young tableaux
(the projective limit of $T_n$ with respect to the natural projections
forgetting the tail of a path). With the topology of coordinatewise convergence
$T$ is a totally disconnected metrizable compact space.
The tail equivalence relation $\sim$ on $T$ is defined as follows:
paths $s=(\mu_1,\mu_2,{\ldots})$ and
$t=(\la_1,\la_2,{\ldots})$ are equivalent
if and only if they are confinal, i.e., 
$\mu_k=\la_k$ for all sufficiently large $k$.
Denote by $[t]_n\in T_n$ the initial segment of length $n$
of a tableau $t\in T$. Given a finite tableau $u\in T_n$, denote
by $C_u=\{t:[t]_n=u\}$ the corresponding cylinder set;
for $\la\in\yn$, let $C_\la=\{t:t_n=\la\}=\cup_{u\in T_\la}C_u$
be the set of all paths passing through $\la$.

A measure $M$ on the space $T$ is called {\it central}
if for any $n\in\N$, $\la\in\yn$ and any tableaux $u,v\in T_\la$ of shape
$\la$, the measures of the corresponding cylinder sets coincide:
$M(C_u)=M(C_v)$ (equivalently, a measure $M$ is central if
it is invariant under any transformations changing the beginnings of tableaux).
One can easily see that the cotransition probabilities
$\frac{M(C_\la\cap C_\La)}{M(C_\La)}$ of a central measure
depend only on the pair of diagrams
$\la\subset\La$ (and do not depend on the measure) and are equal to
$\frac{\dim\la}{\dim\La}$ (the relative dimension of the representation
$\pi_\la$ in $\pi_\La$).
Central measures play a very important role in the representation theory
of the symmetric groups.\footnote{In the theory of
dynamical systems, such measures are called measures of maximum
entropy.}

A measure $M$ on the space $T$ is called {\it Markov}
if for every $n\in\N$ the following condition holds:
for any diagrams $\la\in\yn$ and 
$\La\in{\mathbb Y}_{n+1}$ such that $\La\subset\la$
and for any path $u\in T_\la$, the events
$C_u$ (``the past'') and $C_\La$ (``the future'') are independent
given
$C_\la$ (``the present''). In other words, 
a random tableau $t=(\la_1,\la_2,{\ldots} )$,
regarded as a sequence of random variables $\la_n$, where $\la_n$ 
takes values in the set $\yn$ of Young diagrams with $n$ cells,
is a Markov chain in the ordinary sense. In terms of
transition probabilities, this means that
the transition probability $\frac{M(C_\La\cap C_u)}{M(C_u)}$ 
depends only on the form $\la$ of a tableau $u$, but not on the tableau itself. 
Note that the ``forward'' and ``backward'' Markov properties are
equivalent, so that the definition of a Markov measure can be formulated
in a similar way in terms of cotransition probabilities.

It is easy to see that every central measure is Markov.
Thus central measures form a class of Markov measures with
fixed cotransition probabilities equal to
$\frac{\dim\la}{\dim\La}$.

\section{Markov vectors and Markov representations}
\label{sec:markov}

Consider a cyclic representation $\pi$ of the group 
$\sn$ in a space $V$ that has a simple spectrum (i.e.,
decomposes into the sum of pairwise nonequivalent irreducible representations)
and a unit cyclic vector $\xi$ in this representation.
Consider a pair of diagrams $\mu \supset \la$ with $n$
and $n-1$ cells, respectively, and assume that
$\pi$ contains a subrepresentation equivalent to
$\pi_\mu$ (for convenience, we will denote it by the same symbol).
Project the vector $\xi$ to the space $V_\mu$ of
the representation $\pi_\mu$ and denote the obtained vector by
$\xi_{\mu}$; then project $\xi_{\mu}$ to the space  (contained in $V_\mu$)
of the representation
$\pi_\la$ of the group ${\mathfrak S}_{n-1}$
and denote the obtained vector by $\xi_{\mu,\la}$. 
Let us call the ratio of the squared norms
 \be
 (||\xi_{\mu,\la}||/||\xi_{\mu}||)^2
 \label{cotrans}
 \ee
the cotransition probability of the pair $\mu$, $\la$.

Let us define a measure on the space of all tableaux $t$ with diagrams
corresponding to the representations occurring in the representation $\pi$
of the group $\sn$ as follows:
the probability of a tableau (i.e.,  a path in the Young graph) is equal to the product
of cotransition probabilities along the whole path.
It is easy to see that this measure is well defined and
coincides with the spectral measure of the vector $\xi$
with respect to the Gelfand--Tsetlin algebra; and the values (\ref{cotrans})
are exactly the cotransition probabilities of this spectral measure in the ordinary sense.

Note also that in our considerations the cyclic vector can be
multiplied by a scalar factor of modulus $1$ --- this 
does not change the spectral measure; any other change of the cyclic vector
changes the measure.

\begin{definition}
We will say that a cyclic vector $\xi$ is {\it Markov}
if its spectral measure with respect to the Gelfand--Tsetlin algebra
is Markov.
\end{definition}

Thus a cyclic vector is Markov if for every
$k<n$ and every diagram $\la\in{\mathbb Y}_k$, the probability
of any tableau with this diagram does not depend 
on the continuation of this tableau to the level $n$.
In terms of representations and cyclic vectors, this means
that the norm of the projection of the cyclic vector to 
the subspace of the representation of the group 
${\mathfrak S}_k$ equivalent to $\pi_\la$
does not depend on the way in which we have arrived at this subspace. 

Now we use the following simple lemma from representation theory.

\begin{lemma}
\label{l:cyclic}
Assume that in a finite-dimensional Hilbert space
$H$ there is a unitary representation of a group $G$ that
is primary, i.e., decomposes into the direct
(not necessarily orthogonal) sum 
$H=H_1\oplus H_2\oplus{\ldots} \oplus H_n$ 
of equivalent irreducible representations, and in each of them
there is a cyclic vector
$v_i\in H_i$, $i=1,{\ldots} ,n$.
Then the following two assertions are equivalent:

{\rm 1.} For any $i,j$, there exists an isometry
$T_{i,j}:H_i\to H_j$ intertwining the corresponding representations 
such that $T_{i,j}v_i=v_j$.

{\rm 2.} In the cyclic hull of the vector $v=\sum v_i$,
the representation is irreducible.
\end{lemma}

Applying this lemma to the direct sum of the subrepresentations of $\pi$
equivalent to $\pi_\la$, we see that the norms of the projections
of the cyclic vector to the corresponding subspaces coincide
if and only if $\pi_\la$ has multiplicity $1$ 
in the decomposition
of the representation of ${\mathfrak S}_k$ in the cyclic hull 
${\mathfrak S}_k\xi$ of $\xi$. Since, as shown above, the coincidence
of these norms is in turn equivalent to the Markov property, 
we have proved the following lemma on characterization
of Markov vectors.

\begin{lemma}
\label{l:mark}
Let $\pi$ be a unitary representation of the group
$\sn$ with simple spectrum. A cyclic vector
$\xi$ of the representation $\pi$ is Markov if and only if
for every $k<n$, the representation of the group
${\mathfrak S}_k$ in the cyclic hull
${\mathfrak S}_k\xi$ of $\xi$ with respect to 
${\mathfrak S}_k$ has a simple spectrum.
\end{lemma}

Note that the described procedure of constructing the spectral measure
does not determine a measure on diagrams, because in general the probabilities
of the same diagram in different tableaux do not coincide.
If they do coincide, then the spectral measure is central.

\medskip
Now let us consider representations of the infinite symmetric group $\si$.

If we are given a quasi-invariant measure $\mu$ on the space of Young tableaux
$T$ and a $1$-cocycle $c$ on pairs of confinal paths
taking values in the group of complex numbers of modulus $1$,
then we can construct a unitary representation of the group $\si$
in the space $L^2(T,\mu)$ as follows (see, e.g., \cite{Itogi}). Recall
that the Fourier transform allows one to realize the group algebra
$\C[\si]$ of the infinite symmetric group as the cross product
constructed from the commutative algebra of functions on the space of tableaux $T$
(Gelfand--Tsetlin algebra) and the tail equivalence relation.
The desired representation is given by
\be
L_g h(s)=\sum_{t\sim s}\sqrt{\frac{d\mu(s)}{d\mu(t)}}\hat g(s,t)c(s,t)h(t),
\qquad \phi\in L^2(T,\mu),
\label{repr}
\ee
where $\hat g$ is the function on pairs of confinal paths
corresponding to an element $g\in\si$ 
(the Fourier transform of $g$). Note that the cocycle is trivial on
the space of finite tableaux.

\begin{definition} {\rm
A representation of the infinite symmetric group
$\si$ is called
{\it simple} if it is the inductive\footnote{Recall that the
inductive limit of unitary representations $\pi_k$  
of finite groups $G_k$ forming  an inductive chain
$G_1\subset G_2\subset{\ldots}$  
is the representation of the group
$G=\cup G_k$ in the Hilbert space that is the completion of the chain
of spaces of the representations $\pi_k$ with (equivariant) isometric embeddings.}
limit of representations of the finite symmetric groups
$\sn$ with simple spectrum.}
\end{definition}

\noindent{\bf Remarks. 1.}
Strictly speaking, we should say ``a simple representation
with respect to the approximation of the group
$\si$ by the sequence of finite groups $\sn$, $n=1,2,{\ldots}$,
with standard embeddings.'' Another approximation will give another
class of simple representations. But since in this paper
we use only the standard approximation, we omit this specification 
(see also the remark below on generalized Markov measures).

\medskip
\noindent{\bf 2.} We use the term ``simple representation'' in a much
wider sense than in the paper
\cite{Itogi}, where a representation was called simple if it is the 
inductive limit
of irreducible representations; we will call such representations
{\it elementary}. Of course, elementary representations are simple
in our sense.
\medskip

Note that a representation with simple spectrum is cyclic, i.e., contains
a cyclic vector.

\begin{definition}
{\rm  A representation $\pi$ of the group
$\si$ with simple spectrum  is called
{\it Markov} if the space of $\pi$ contains a cyclic vector
whose spectral measure (with respect to the Gelfand--Tsetlin algebra)
is Markov.
Note that a representation with simple spectrum is Markov if and only if
the measure $\mu$ in its realization~(\ref{repr}) is Markov.}
\end{definition}

\begin{theorem}
\label{th:mark}
A representation of the infinite symmetric group
is Markov if and only if it is simple.
\end{theorem}

\begin{proof}
If $\xi$ is a Markov cyclic vector, then for every $n$ the representation
of the group $\sn$ in the cyclic hull $\sn\xi$ is Markov and hence, by 
Lemma~\ref{l:mark}, has a simple spectrum. Thus the representation
is the inductive limit of representations with simple spectrum.
Conversely, if we have the inductive limit of representations with simple spectrum, 
then, by Lemma~\ref{l:mark}, the vector obtained by the successive embeddings
from the unit vector in the original one-dimensional representation
of ${\mathfrak S}_1$ is a Markov cyclic vector.
\end{proof}

Consider a Markov representation $\pi$ and a Markov cyclic vector
$\xi$. The spectrum of the representation of $\sn$ in its cyclic hull
$\sn\xi$ is simple (by Lemma~\ref{l:mark}), and we obtain an approximation of
$\pi$ by representations with simple spectrum, but the action of 
$\sn$ in this representation may differ from the standard
Young form by a factor of modulus $1$.
Of course, correcting the basis by such a factor 
(introducing the ``phase''), we can obtain the standard action, 
but it may happen that there is no convergent system of factors
and, consequently, the limit action involves a cocycle
$c$, as in~(\ref{repr}). The question of convergence of factors
is precisely the question of 
whether the cocycle is cohomological to the identity one.
It is well known that there are many cocycles that are not
cohomological to the identity one. Thus, for a given Markov measure,
there are many nonequivalent representations differing by a cocycle.
The cocycle ``measures the deviation'' of the given realization 
of a representation from its standard realization. We do not dwell
on this interesting questions, because in all realizations considered below
(tensor model, induced representations), the cocycle is trivial.

It is not difficult to see how Theorem~\ref{th:mark} can be extended  
to the case of generalized Markov measures. Assume that the spectral 
measure of the representation is a generalized Markov measure with Markov
intervals $N_1,N_2,\ldots$ 
(i.e., for every $k$, if we fix the diagrams with
$N_{k-1},{\ldots},N_k-1$ cells, then the diagrams preceding
$N_{k-1}$ and those following $N_k-1$ are independent); then
the corresponding representation is simple, but with respect to
a sparse chain of subgroups
${\mathfrak S}_1,{\mathfrak S}_{N_1},
{\mathfrak S}_{N_2},{\ldots}$. In other words, it is the limit
of representations with simple spectrum of the sequence of groups
${\mathfrak S}_{N_k}$, $k \to \infty$. As mentioned in the introduction,
every quasi-invariant measure on the space of tableaux is equivalent
to a generalized Markov measure --- the proof of this fact is the same
as for a similar assertion for quasi-invariant measures in the space 
of one-sided sequences (see \cite{V}). Therefore, {\it each irreducible
representation of the group $\si$ with simple spectrum
(with respect to the Gelfand--Tsetlin algebra) is the limit
of finite-dimensional representations with simple spectrum
of a certain sequence of groups ${\mathfrak S}_{N_k}$, $k=1,2,{\ldots} $}.

\section{The tensor model of two-row representations}
\label{sec:tensor}

In this section, we describe 
the so-called tensor model of two-row representations of the symmetric groups,
which was suggested by the first author and studied in \cite{Nik}
(see also \cite{TsV}).

For $0\le k\le n$, denote by 
$\Fnk$
the set of $k$-element subsets in
$\{1,{\ldots} ,n\}$.
Given $I=\{i_1,{\ldots} , i_k\}\in\Fnk$, let
$x_I=x_{i_1}\cdot\ldots\cdot x_{i_k}$.

Denote by $A_{n,k}=\{\sum_{I\in\Fnk} c_I x_I\}$
the vector space of square-free homogeneous forms of degree
$k$ in $n$ variables. 
This space can also be identified with 
the space of symmetric tensors of rank $k$
with zero diagonal components over the $n$-dimensional space
(a form
$f=\sum_{I\in\Fnk} c_I x_I$ is identified with the tensor
$\{T_{j_1,{\ldots} ,j_k}\}_{j_1,{\ldots} ,j_k=1}^n$, where 
$T_{j_1,{\ldots} ,j_k}
=c_{\{j_1,{\ldots} ,j_k\}}$ if the indices 
$j_1,{\ldots} ,j_k$ are pairwise distinct and $T_{j_1,{\ldots} ,j_k}=0$
otherwise).

Denote by $\|\cdot\|$ the standard scalar product in the space
of forms (tensors) given by
\be
\|f\|^2=\sum_{I\in\Fnk} |c_I|^2,\qquad
f=\sum_{I\in\Fnk} c_I x_I\in A_{n,k}.
\label{scal}
\ee

Let
$A^0_{n,k}$ be the subspace of $A_{n,k}$ defined as
$$
A^0_{n,k}=\left\{\sum c_I x_I \in A_{n,k}\,|\, \sum_{j\not\in J} c_{J\cup j}
 = 0\text{ for every }J\in F_{n,k-1} \right\}.
$$

There is a natural action of the symmetric group $\sn$ on the space
 $A_{n,k}$ by substitutions of indices: given $\sigma\in \sn$,  
 $$
 \sigma\cdot\sum_{I\in\Fnk}c_Ix_I=\sum_{I\in\Fnk}c_{I}x_{\sigma I},\quad\mbox{where}\quad \sigma\{i_1,{\ldots} , i_k\}
=\{\sigma(i_1),{\ldots} , \sigma(i_k)\},
 $$
or, in tensor form, $\sigma\{T_{j_1,{\ldots} ,j_k}\}=\{T'_{j_1,{\ldots} ,j_k}\}$, where
$T'_{j_1,{\ldots} ,j_k}=T_{\sigma^{-1}(j_1),{\ldots} ,\sigma^{-1}(j_k)}$.
It is easy to see that the subspace $A^0_{n,k}$ is invariant
under this action. Note that the spaces $A^0_{n,k}$ and $A^0_{n,n-k}$
(as well as $A_{n,k}$ and $A_{n,n-k}$) are naturally isometric
and the corresponding representations of $\sn$
are equivalent.

Realizations of representations of the symmetric groups in spaces 
of tensors (square-free homogeneous forms) are called
{\it tensor realizations}. The following theorem 
is a well-known fact, which in terms of the tensor model was proved in \cite{Nik}.

\begin{theorem} 
Let $k\le n/2$. 

{\rm(1)} 
The representation of the symmetric group
$\sn$ in the space $A^0_{n,k}$ 
(and in the space $A^0_{n,n-k}$) is equivalent to
the irreducible representation $\pi_{n-k,k}$
corresponding to the two-row diagram $\la_{n,k}=(n-k,k)$ with rows 
of lengths $n-k$ and $k$.

{\rm(2)} The representation of the symmetric group
$\sn$ in the space $A_{n,k}$ is 
equivalent to the multiplicity-free direct sum of $\pi_{n-l,l}$ over all
$l=0,1,{\ldots} ,k$. In particular,
the representation of $\sn$ in $A_{n,[n/2]}$ is equivalent to the
multiplicity-free direct sum of all two-row representations:
$$
A_{n,[n/2]}\simeq\bigoplus_{k=0}^{[n/2]}\pi_{n-k,k}.
$$
\label{th:nik}
\end{theorem}

\section{The Gelfand--Tsetlin basis in tensor realizations}
\label{sec:gz}

In this section, we find explicit formulas for the Gelfand--Tsetlin basis
in the tensor realizations of two-row representations of the symmetric groups.

First note that the space $A_{n,k}^0$ can be also defined as follows.

\begin{lemma}
\label{l:0}
The space $A_{n,k}^0$ is the subspace of $A_{n,k}$ consisting
of forms that are invariant under simultaneous shifts of all variables
by a constant. It is spanned by functions of the form
$(x_{i_1}-x_{j_1}){\ldots}(x_{i_k}-x_{j_k})$,
where all indices $i_1,{\ldots} ,i_k,j_1,{\ldots} ,j_k$
are pairwise distinct.
\end{lemma}

Of course, functions of the form
$(x_{i_1}-x_{j_i}){\ldots}(x_{i_k}-x_{j_k})$ are
linearly dependent, so that they
form an overfull system in $A_{n,k}^0$. We will call them 
{\it ``pseudo-monomials.''}

Note that a two-row Young tableau $u$ is uniquely determined by the sequence
$p_1<p_2<{\ldots}<p_k$ of the elements of its second row.

\begin{theorem}
\label{th:gz0}
Let $u\in T_n$ be a two-row Young tableaux, and denote by 
$p_1<p_2<{\ldots}<p_k$
the elements of its second row. Then the element $h^0_u$ of
the Gelfand--Tsetlin basis in $A^0_{n,k}$
corresponding to the tableau $u$ 
is given by the formula
\begin{equation}
\label{ht0}
h_u^0=c_u^0\sum_{i_1,{\ldots} ,i_k}
(x_{i_1}-x_{p_1}){\ldots}(x_{i_k}-x_{p_k}),
\end{equation}
where $c^0_u$ is a normalizing
constant and
the sum runs over all indices $i_1,{\ldots}, i_k$ from $1$ to $n$ such that
$i_j<p_j$ for all $j=1,{\ldots} ,k$ and all indices $i_1,{\ldots} ,i_k,p_1,{\ldots} ,p_k$ 
are pairwise distinct.
\end{theorem}

\begin{proof}
It suffices to check that for every  $l=1,{\ldots} ,n$, the form
$h_u^0$ is an eigenfunction for the action of 
the YJM-element\footnote{Recall that the $l$th Young--Jucys--Murphy (YJM) element
is defined as $X_l=(1,l)+(2,l)+{\ldots} +(l-1,l)\in{\mathbb C}[\sn]$, 
where, as usual, $(i,j)$
stands for the transposition permuting the elements $i$ and $j$. The
YJM-elements form a multiplicative basis of the Gelfand--Tsetlin
algebra $\gzn$. Concerning the Gelfand--Tsetlin
algebras, YJM-elements, and their role in the representation theory
of the symmetric groups, see \cite{VO}.}
$X_l$ with
eigenvalue equal to $c_l(u)$, the contents of the cell of $u$ containing $l$.
This can be done by direct combinatorial calculations.
\end{proof}

One can prove that the normalizing constant in (\ref{ht0}) equals 
\be
\label{ct0}
c^0_u=\frac{1}{\bigl(\prod_{j=1}^k(p_j-2j+1)(p_j-2j+2)\bigr)^{1/2}}.
\ee

\smallskip\noindent{\bf Example 1.} Consider the tableau $u$ of shape 
$\la=\la_{n-k,k}$ in which 
the second row contains the numbers $2,4,{\ldots} ,2k$. Then 
$$
h_u^0=c_u^0\cdot (x_1-x_2)(x_3-x_4){\ldots} (x_{2k-1}-x_{2k}),
\qquad (c_u^0)^2=\frac1{2^k}.
$$
This is the only tableau of shape $\la$
such that the corresponding Gelfand--Tsetlin
element $h_u^0$ is a pseudo-monomial; 
we will call it the {\it good} tableau of shape $\la$.
For all other tableaux of the 
same shape, $h_u^0$ 
is a linear combination of pseudo-monomials. 
\medskip

Denote by $H_{n,m}^k$ the subspace in $A_{n,m}$ corresponding to
the representation $\pi_{{n-k,k}}$. Thus
\be
A_{n,m}=\bigoplus_{k=0}^{m} H_{n,m}^k.
\label{dec}
\ee
Let $\psi_n^l$ be the linear operator that acts on monomials as follows:
$$
\psi_n^lx_I=x_I\sum_{j_1,{\ldots} ,j_l\notin I\atop\mbox{\scriptsize distinct}}x_{j_1}{\ldots} x_{j_l}.
$$
As follows from the results of \cite{Nik}, 
\be
H_{n,m}^k=\psi_n^{m-k}A_{n,k}^0,
\label{psi}
\ee
and $\psi_n^{m-k}$ is an isomorphism of $H_{n,m}^k$ and $A_{n,k}^0$
intertwining the corresponding representations.

Using (\ref{psi}), we can obtain 
the following assertions from Lemma~\ref{l:0} and Theorem~\ref{th:gz0}.

\begin{lemma}
The space $H_{n,m}^k$ is spanned by functions of the form 
$(x_{i_1}-x_{j_1}){\ldots}(x_{i_k}-x_{j_k})x_{s_1}{\ldots}x_{s_{m-k}}$,
where all indices $i_1,{\ldots} ,j_k,j_1,{\ldots} ,j_k, s_1,{\ldots} ,s_{m-k}$
are pairwise distinct.
\end{lemma}

\begin{theorem}
Let $u\in T_n$ be a two-row Young tableaux, and denote by $p_1<p_2<{\ldots}<p_k$
the elements of its second row. Then the element $h_u$
of the Gelfand--Tsetlin basis in $H_{n,m}^k$
corresponding to the tableau $u$ is given by the formula
\begin{equation}
\label{ht}
h_u=c_u\sum
(x_{i_1}-x_{p_1}){\ldots}(x_{i_k}-x_{p_k})
x_{s_1}{\ldots}x_{s_{m-k}},
\end{equation}
where $c_u$ is a normalizing constant and
the sum runs over all indices $i_1,{\ldots}, i_k$,
$s_1,{\ldots} ,s_{m-k}$ from $1$ to $n$ such that
$i_j<p_j$ for all $j=1,{\ldots} ,k$ and all indices
$i_1,{\ldots} ,i_k,p_1,{\ldots} ,p_k,s_1,{\ldots} ,s_{m-k}$ 
are pairwise distinct.
\end{theorem}

\smallskip\noindent{\bf Example 2.} For the good tableau $u$ of shape 
$\la=\la_{n-k,k}$ 
considered in Example~1, we have
\be
\label{pure}
h_u=c_u\cdot (x_1-x_2)(x_3-x_4){\ldots} (x_{2k-1}-x_{2k})
\sum_{s_1,{\ldots} ,s_{m-k}\in[2k+1..n]\atop\mbox{\scriptsize distinct}}
x_{s_1}{\ldots} x_{s_{m-k}},\qquad c_u^2=\frac1{2^kC_{n-2k}^{m-k}}.
\ee

\section{Induced representations}
\label{sec:induced}

In this section, we consider the class of irreducible
representations of the infinite symmetric group $\si$ induced from two-block
Young subgroups. Namely, consider a partition $\mathbb N=A\cup B$
of the set of
positive integers into two subsets 
and the representation of $\si$ induced from the identity representation
of the subgroup ${\mathfrak S}_{A}\times
{\mathfrak S}_{B}$. It is not difficult to prove that this 
representation is irreducible. It
can also be described in terms of the tensor model as follows.

A partition $\mathbb N=A\cup B$ of the set of positive integers into two subsets 
is uniquely determined  by 
an infinite sequence 
$\xi=\xi_1\xi_2{\ldots}$ of $0$'s and
$1$'s (an ``infinite tensor''), 
where $\xi_i=1$ if $i\in A$, and $\xi_i=0$ if $i\in B$.  
Then the induced representation in question is equivalent to
the natural substitutional representation of $\si$ 
on infinite sequences in the cyclic hull
of the sequence $\xi$, which we will denote
by $\pi_\xi$.
Note that the orbit of $\xi$ is the discrete set $O_\xi$ of infinite sequences
of $0$'s and $1$'s eventually coinciding with $\xi$,
and $\pi_\xi$ is a unitary representation of $\si$ 
in the space $l^2(O_\xi)$.

For simplicity, it is convenient to assume that the number of 
$1$'s among the first $n$ elements of $\xi$ does not exceed $n/2$.
It is not difficult to see that an arbitrary case can be reduced
to this one, but we omit the corresponding technical details.

Consider 
the cyclic hull $\sn\xi$ of the sequence $\xi$ with respect to 
the finite symmetric group $\sn$.
It can be naturally identified with the space $A_{n,m}$, 
where $m=m(n)$ is the number of $1$'s among the first
$n$ elements of $\xi$. Moreover, the scalar product induced from
$l^2(O_\xi)$ coincides with the standard scalar product~(\ref{scal}) 
in $A_{n,m}$. Thus the representation $\pi_\xi^n$ of the group $\sn$ in $\sn\xi$ 
is unitarily equivalent to the tensor representation in the space
$A_{n,m}$, the cyclic vector $\xi$ corresponding to the monomial
$x_{i_1}{\ldots} x_{i_m}\in A_{n,m}$, where
$i_1,{\ldots} ,i_m$ are the numbers of positions from $1$ to $n$
at which $\xi$ has $1$'s.
We also have a natural embedding 
$\iota_n:\pi_\xi^n\hookrightarrow\pi_\xi^{n+1}$ defined as follows: if $\xi_{n+1}=0$, then
$\iota_n$ is the identical embedding $A_{n,m}\hookrightarrow A_{n+1,m}$;
and if $\xi_{n+1}=1$, then $\iota_n:A_{n,m}\to A_{n+1,m+1}$ is the multiplication by $x_{n+1}$,
i.e., $\iota_n f=x_{n+1}f$.
The following lemma is obvious.

\begin{lemma}
The induced representation $\pi_\xi$ of the group $\si$
is the inductive limit of
the tensor representations $\pi_\xi^n$ of the groups $\sn$.
\end{lemma}

The main result of this section is the following Theorem~\ref{th:spectral}. Note that 
it easily follows from the tensor realization, 
Theorem~\ref{th:nik}, and Lemma~\ref{l:mark}
that the representation $\pi_\xi$ is Markov, so that the main 
part of the theorem is the computation of the transition
probabilities of the corresponding spectral measure.

\begin{theorem}
\label{th:spectral}
The spectral measure $\mu_\xi$
of the cyclic vector $\xi$ in the representation $\pi_\xi$
with respect to the Gelfand--Tsetlin algebra is a Markov measure on the
space of infinite Young tableaux $T$, and its transition probabilities 
are given by the 
following formula. Denote by $m(n)$ the number of $1$'s among the first
$n$ elements of $\xi$.

If $\xi_{n+1}=0$, then
\be
\label{8}
\Prob(\la_{n,k},\la_{n+1,k})=\frac{n-m(n)-k+1}{n-2k+1},\qquad
\Prob(\la_{n,k},\la_{n+1,k+1})=\frac{m(n)-k}{n-2k+1}.
\ee

If $\xi_{n+1}=1$, then
\be
\label{9}
\Prob(\la_{n,k},\la_{n+1,k})=\frac{m(n)-k+1}{n-2k+1},\qquad
\Prob(\la_{n,k},\la_{n+1,k+1})=\frac{n-m(n)-k}{n-2k+1}.
\ee
\end{theorem}

\begin{proof}
As mentioned above, for all $n$ the representation $\pi_\xi^n$ is realized in the space
$A_{n,m}$, where $m=m(n)$ is the number of $1$'s among the first
$n$ elements of $\xi$. 
In particular, for all $l\le n$, 
the representation of ${\mathfrak S}_l$ in the 
corresponding cyclic hull ${\mathfrak S}_l\xi$
of $\xi$ has a simple spectrum by Theorem~\ref{th:nik}. It follows from
Lemma~\ref{l:mark} that $\mu_\xi$ is a Markov measure.

Let us find explicit formulas for 
the transition probabilities of $\mu_\xi$.
As follows from the branching rule for irreducible representations
of the symmetric groups,
for each vector $f\in H_{n,m}^k$, we have
$\iota_nf=f_{n+1,k}+f_{n+1,k+1}$, where
$f_{n+1,k}\in H_{n+1,m'}^k$, $f_{n+1,k+1}\in H_{n+1,m'}^{k+1}$, 
and $m'=m(n+1)=m+\xi_{n+1}$.

\begin{lemma}
\label{l:dec}
Let $f\in H_{n,m}^k$. By {\rm(\ref{psi})}, $f=\psi_n^{m-k}f_0$,
where $f_0\in A_{n,k}^0$. Then

$\bullet$ if $\xi_{n+1}=0$,
\begin{eqnarray*}
f_{n+1,k}&=&\frac{n-m-k+1}{n-2k+1}\bigl(f+x_{n+1}\psi_n^{m-k-1}f_0\bigr),\\
f_{n+1,k+1}&=&\frac{1}{n-2k+1}\bigl((m-k)f-(n-m-k+1)x_{n+1}\psi_n^{m-k-1}f_0\bigr);
\end{eqnarray*}

$\bullet$ if $\xi_{n+1}=1$,
\begin{eqnarray*}
f_{n+1,k}&=&\frac{m-k+1}{n-2k+1}\bigl(x_{n+1}f+\psi_n^{m-k+1}f_0\bigr),\\
f_{n+1,k+1}&=&\frac{1}{n-2k+1}\bigl((n-m-k)x_{n+1}f-(m-k+1)\psi_n^{m-k+1}f_0\bigr).
\end{eqnarray*}
\end{lemma}
\begin{proof}
One can easily check that the forms $f_{n+1,k}$ and $f_{n+1,k+1}$ lie
in the appropriate spaces. Obviously, their sum equals $\iota_nf$.
\end{proof}

\begin{corollary}
\label{cor:norms}
Let $f\in H_{n,m}^k$. Then

$\bullet$ if $\xi_{n+1}=0$, then
\begin{eqnarray*}
\|f_{n+1,k}\|^2=\frac{n-m-k+1}{n-2k+1}\|f\|^2,\qquad
\|f_{n+1,k+1}\|^2=\frac{m-k}{n-2k+1}\|f\|^2.
\end{eqnarray*}

$\bullet$ if $\xi_{n+1}=1$, then
\begin{eqnarray*}
\|f_{n+1,k}\|^2=\frac{m-k+1}{n-2k+1}\|f\|^2,\qquad
\|f_{n+1,k+1}\|^2=\frac{n-m-k}{n-2k+1}\|f\|^2.
\end{eqnarray*}
\end{corollary}
\begin{proof}
The desired formulas follow from Lemma~\ref{l:dec} by 
straightforward calculations taking into account that
$(f_{n+1,k},f_{n+1,k+1})=0$; $\|f_{n+1,k}\|^2+\|f_{n+1,k+1}\|^2=\|f\|^2$;
$(f,x_{n+1}\psi_n^{m-k-1}f_0)=0$ in the case $\xi_{n+1}=0$, and
$(x_{n+1}f,\psi_n^{m-k+1}f_0)=0$ in the case $\xi_{n+1}=1$
(the first two relations are valid
by definitions, and the last one holds because one of the vectors consists of monomials
containing $x_{n+1}$ and the other one consists of monomials not
containing $x_{n+1}$).
\end{proof}

As follows from the proof of Lemma~\ref{l:mark}, 
$\Prob(\mu\,|\,\la)=\bigl(\frac{\|f_\la\|}{\|f_\mu\|}\bigr)^2$, so that
formulas~(\ref{8}), (\ref{9}) follow from Corollary~\ref{cor:norms}.
Theorem~\ref{th:spectral} is proved.
\end{proof}

Note that all spectral measures of induced representations 
considered in Theorem~\ref{th:spectral}
are not central (except for the trivial case when one of the sets
in the partition is empty).

\begin{corollary} The operator $\psi_n^{m-k}:A^0_{n,k}\mapsto H_{n,m}^k$
is an isometry up to a constant. Namely, 
for $f_0\in A^0_{n,k}$, 
\be
\label{psinorm}
\|\psi_n^{m-k}f_0\|^2=C_{n-2k}^{m-k}\cdot\|f_0\|^2.
\ee
\end{corollary}
\begin{proof}
Follows by simple calculations based on Lemma~\ref{l:dec} and Corollary~\ref{cor:norms}.
\end{proof}

In particular, we obtain from~(\ref{ct0}), (\ref{psi}), and (\ref{psinorm})
that the coefficients $c_u$ in~(\ref{ht}) satisfy
$$
\frac{1}{c_u^2}=C_{n-2k}^{m-k}\cdot\prod_{j=1}^k(p_j-2j+1)(p_j-2j+2).
$$

\smallskip\noindent{\bf Example 3.}
Let $\xi=0101{\ldots}$. Then $m(n)=[n/2]$, and the formulas for transition 
probabilities take the following form:

$\bullet$ if $n$ is odd,
\be
\Prob(\la_{n,k},\la_{n+1,k})=\frac{n-2k+2}{2(n-2k+1)},\qquad
\Prob(\la_{n,k},\la_{n+1,k+1})=\frac{n-2k}{2(n-2k+1)};
\label{12odd}
\ee

$\bullet$ if $n$ is even,
\be
\Prob(\la_{n,k},\la_{n+1,k})=
\Prob(\la_{n,k},\la_{n+1,k+1})=\frac12.
\label{12even}
\ee
These formulas can also be obtained in another way.
Since $\mu_\xi$ is a Markov measure, it suffices 
to find the transition probabilities for one
tableau $u$ of each shape $\la=\la_{n-k,k}$. Let $u$ be
the good tableau of shape $\la$ (see Examples~1 and~2). Note that the
tableau $v$ obtained from a good tableau $u$ by adding the element
$n+1$ to the first row is also good.
Thus, using formula~(\ref{pure}) and observing that in our case  
$\xi$ is the monomial $x_2{\ldots} x_{2n}$, we see that
this monomial occurs exactly once in the right-hand side of~(\ref{pure}),
so that $(\xi,h_u)=c_u$. Hence 
$$
\Prob(\la_{n,k},\la_{n+1,k})=\frac{(\xi,h_v)^2}{(\xi,h_u)^2}=
\frac{C_{n-2k}^{m(n)-k}}{C_{n+1-2k}^{m(n+1)-k}},
$$
which implies the required formula for $\Prob(\la_{n,k},\la_{n+1,k})$, 
taking into account that
$m(n)=[n/2]$. Obviously, 
$\Prob(\la_{n,k},\la_{n+1,k+1})=1-\Prob(\la_{n,k},\la_{n+1,k})$.

\smallskip
It is interesting to compare formulas~(\ref{12odd}), (\ref{12even}) with
the transition probabilities of the ergodic 
{\it central} measure $\mu_\al$
on $T$ corresponding to the Thoma parameters $\al=(1/2,1/2,0,{\ldots})$,
$\beta=0$ (see, e.g., \cite{GB}). 

\begin{lemma}
The transition probabilities of the central measure
$\mu_\al$ are given by~{\rm(\ref{12odd})} for all $n$.
\end{lemma}
\begin{proof}
Let $u\in T_\la$, $\la\in\yn$. Then
$$
\mu_\al(C_u)=s_\la(1/2,1/2)=\frac{1}{2^n}s_\la(1,1)=
\prod_{\square\in\la}\frac{2+c(\square)}{h(\square)},
$$
where $s_\la$ is a Schur function,
$c(\square)$ and
$h(\square)$ are the contents and the hook length of a square
$\square\in\la$, respectively, and we have used the well-known formula for
$s_\la(1,{\ldots} ,1)$
(see \cite[Example I.3.4]{Mac}). The lemma follows by elementary
calculations.
\end{proof}

It is convenient to rewrite formulas~(\ref{12odd}), (\ref{12even})
introducing the change of indices $j=n-2k$. In these terms,
a Young tableau is determined by a sequence 
$(j_1,j_2,{\ldots})$, where  $j_n$ takes the values
$0,1,{\ldots} ,n$,
and the transition probabilities of the measure $\mu_\xi$ are equal to
$$
\Prob(j,j+1)=\frac{j+2}{2(j+1)},\qquad 
\Prob(j,j-1)=\frac{j}{2(j+1)}
$$ 
at an odd moment of time; and
$$
\Prob(j,j+1)=\Prob(j,j-1)=\frac12
$$ 
at an even moment of time.
We see that a random Young tableau governed by the measure $\mu_\xi$
is a trajectory of 
a nonhomogeneous (neither in time nor in space)
random walk on ${\mathbb Z}_+$.
Thus the induced representations of the infinite symmetric group
considered in this paper act in spaces of functions over
trajectories of natural random walks. Explicit formulas for
this action are given by Young's orthogonal form (see, e.g., \cite{JK}).

\end{document}